\documentclass[11pt]{article}

\usepackage[
	a4paper,
	margin=1in
]{geometry}
\usepackage{setspace}
\usepackage{sectsty} 
\usepackage{titlesec}
\usepackage{appendix}

\usepackage{amsmath, amssymb, amsthm}
\usepackage{mathtools}

\usepackage[
	setpagesize=false,
	colorlinks=true,
	linkcolor=BrickRed,
        citecolor=OliveGreen,
        urlcolor=black,
	pdfencoding=auto,
	psdextra,
]{hyperref}
\usepackage{cleveref}

\usepackage{etoolbox} 
\usepackage[shortlabels]{enumitem}
\usepackage{xparse} 

\usepackage{graphicx}
\usepackage[dvipsnames]{xcolor}
\usepackage{tikz}

\usepackage{todonotes}
\usepackage[
    deletedmarkup=sout,
    commentmarkup=uwave,
    authormarkuptext=name
]{changes}

\usepackage{comment}

\setlist[enumerate,1]{label=(\arabic*), ref=(\arabic*)}
\setlist[enumerate,3]{label=(\roman*), ref=(\roman*)}

\allsectionsfont{\boldmath} 

\titlelabel{\thetitle.\quad}

\theoremstyle{plain}
\newtheorem{theorem}{Theorem}[section]
\newtheorem{lemma}[theorem]{Lemma}

\newtheorem{conjecture}[theorem]{Conjecture}
\newtheorem{problem}[theorem]{Problem}

\newtheorem*{claim*}{Claim}

\makeatletter
\newenvironment{claimproof}[1][Proof]{\par
	\pushQED{\qed}%
	
	\normalfont \topsep6\p@\@plus6\p@\relax
	\trivlist
	\item[\hskip\labelsep
	\textit{#1}\@addpunct{.}~]\ignorespaces
}{%
	\popQED\endtrivlist\@endpefalse
}
\makeatother

\theoremstyle{definition}

\newtheorem*{definition*}{Definition}

\newcommand{\ZZ}{\mathbb{Z}}

\newcommand{\calL}{\mathcal{L}}

\newcommand{\calS}{\mathcal{S}}

\newcommand{\ve}{\varepsilon}

\newcommand{\eps}{\varepsilon}

\newcommand{\defeq}{\coloneqq}

\usetikzlibrary{matrix,arrows,decorations.pathmorphing}
\usetikzlibrary{positioning}
\usetikzlibrary{graphs} 
\usetikzlibrary{graphs.standard}

\newcounter{propcounter}
\stepcounter{propcounter}

\title{On a Ramsey--Tur\'{a}n variant of Roth's theorem}

\author{Matija Buci\'{c}\thanks{Department of Mathematics, Princeton University, Princeton, USA. Research supported in part by NSF Award DMS-2349013. Email: \href{mailto:mb5225@princeton.edu}{\nolinkurl{mb5225@princeton.edu}}.} \and 
Micha Christoph\thanks{Department of Mathematics, ETH Z\"{u}rich, Switzerland. Research supported by SNSF Ambizione Grant No. 216071. Email: \href{mailto:micha.christoph@math.ethz.ch}{\nolinkurl{micha.christoph@math.ethz.ch}}.} \and 
Jaehoon Kim\thanks{Department of Mathematical Sciences, KAIST, South Korea.
        \emph{E-mails:} \textbf{\{jaehoon.kim, hyunwoo.lee\}@kaist.ac.kr}. Supported by the National Research Foundation of Korea (NRF) grant funded by the Korea government(MSIT) No. RS-2023-00210430.} \and
Hyunwoo Lee\footnotemark[3]~\thanks{Extremal Combinatorics and Probability Group (ECOPRO), Institute for Basic Science (IBS). Supported by  the Institute for Basic Science (IBS-R029-C4).} \and
Varun Sivashankar\thanks{Department of Mathematics, Princeton University, Princeton, USA. Email: \href{mailto:varunsiva@princeton.edu}{\nolinkurl{varunsiva@princeton.edu}}.} }

\usepackage[square,sort,comma,numbers]{natbib}
\begin{document}
\maketitle

\begin{abstract}
    A classical theorem of Roth states that the maximum size of a solution-free set of a homogeneous linear equation $\mathcal{L}$ in $\mathbb{F}_p$ is $o(p)$ if and only if the sum of the coefficients of $\mathcal{L}$ is $0$. In this paper, we prove a Ramsey--Tur\'{a}n variant of Roth's theorem, with respect to a natural notion of ``structured'' sets introduced by Erd\H{o}s and S\'ark\"ozy in the 1970's. 
    Namely, we show that the following statements are equivalent:
    \begin{enumerate}
        \item[$(a)$] Every solution-free set $A$ of $\mathcal{L}$ in $\mathbb{F}_p$ with $\alpha(\mathrm{Cay}_{\mathbb{F}_p}(A)) = o(p)$ has size $o(p)$. 

        \item[$(b)$] There exists a non-empty \emph{subset} of coefficients of $\mathcal{L}$ with zero sum.
    \end{enumerate}
\end{abstract}


\section{Introduction}\label{sec:intro}

In the last decades, many surprising connections between additive combinatorics and extremal (hyper)graph theory have been discovered. A classical example is Schur's theorem~\cite{Schur1917}, which states that every coloring of $[n],$ for sufficiently large $n,$ using $r$ colors contains a monochromatic solution to the equation $x + y = z$. The proof of Schur's theorem relies on the finiteness of the (multicolor) Ramsey number of a triangle, and the quantitative bounds for Schur's theorem are closely related to the quantitative bounds for the Ramsey numbers. 
This tantalizing connection permeates the area of Arithmetic Ramsey theory. 
One of the most notable theorems in early Arithmetic Ramsey theory is Van der Waerden's theorem~\cite{van1927beweis}. It states that every coloring of $[n]$, for sufficiently large $n$, using $r$ colors contains a monochromatic $k$-term arithmetic progression. 

Inspired by Van der Waerden's theorem, Erd\H{o}s and Tur\'{a}n~\cite{Erdos-Turan} conjectured the following: For every $\ve > 0$ and $k \geq 3$, there exists $n_0 = n_0(k, \ve)$ such that for all $n \geq n_0$, every subset $A\subseteq [n]$ with $|A| \geq \ve n$ contains a $k$-term arithmetic progression. The first nontrivial case of this conjecture, the case $k = 3$, was proved by Roth~\cite{Roth-I}, using Fourier analytic methods. Subsequently, Szemer\'{e}di~\cite{Szemeredi} famously found a graph-theoretic proof of the full conjecture, nowadays referred to as Szemer\'{e}di's theorem. The proof led to the development of the regularity lemma, one of the most powerful tools in modern extremal graph theory.

Recently, besides analogues of Ramsey or Tur\'{a}n-type problems in additive combinatorics, there has been extensive work in additive combinatorics inspired by topics from extremal graph theory, such as common and Sidorenko linear equations~\cite{Saad-Wolf,Sidorenko-equation}, forcing and norming linear equations~\cite{norming-equations}, to name just a few. We recommend a wonderful book~\cite{Zhao-book} by Zhao, for more details, and further connections between extremal graph theory and additive combinatorics. 

In this article, we investigate a new such direction, namely a Ramsey--Tur\'{a}n variant of Roth's theorem. 
In the classical Tur\'an problem, one asks for the maximum number of edges in a graph $G$ having no subgraph $F$. If $F$ has chromatic number at least three, then the asymptotically tight lower bound is attained by a balanced blow-up of a complete graph of order $\chi(F)-1$. This extremal example is highly structured, leading to the very natural question of whether this is necessary, and in particular whether one may obtain better bounds for more ``typical'' $F$-free host graphs $G$. A classical way of ensuring typicality, introduced by Erd\'{o}s and S\'{o}s~\cite{Erdos-Sos} in 1970, is to restrict attention only to host graphs $G$ with sublinear independence number. This led to the development of the extensively studied Ramsey--Tur\'{a}n theory, (see e.g.\ a survey~\cite{Simonovits-Sos}). 

In this paper, we study a Ramsey--Tur\'{a}n variant of the following classical theorem of Roth~\cite{Roth-2} classifying density-regularity of homogeneous linear equations over $\mathbb{F}_p$. 
We say $A \subseteq \mathbb{F}_p$ is a \emph{solution-free} set of $\calL: \sum_{i\in [k]} c_i x_i = 0$ if there is no $k$-tuple $(a_1, \dots, a_k)\in A^k$ of distinct\footnote{We note here, that we follow Roth in requiring the solution to use distinct elements of $A$, although most of the results remain the same and equally interesting if one drops this assumption.} elements of $A$ satisfying the equation, namely $\sum_{i\in [k]} c_i a_i = 0$.

\begin{theorem}[Roth~\cite{Roth-2}]\label{thm:roth}
    Let $\calL : \sum_{i\in [k]} c_i x_i = 0$ be a homogeneous linear equation with $k \ge 3$ and $ c_1,\ldots, c_k \in \mathbb{Z} \setminus \{0\}.$ Let further $p$ be a prime. Then, the following are equivalent
    \begin{enumerate}
        \item[$(a)$] Every solution-free set $A$ of $\calL$ in $\mathbb{F}_p$ has size $o(p)$.
        \item[$(b)$] $\sum_{i\in [k]} c_i = 0$. 
    \end{enumerate}
\end{theorem}

Similarly, as in the case of Tur\'an's theorem, the standard examples of solution-free sets establishing the only if part of Roth's theorem are highly structured. A natural way of avoiding structured examples, in particular with the ideas from Ramsey--Tur\'an theory in mind, is to only consider solution-free sets with small independence number in their Cayley graphs. We note that this notion of avoiding structure (stated in a different language) dates back to the works of Erd\H{o}s and S\'ark\"ozy \cite{erdHos1977differences} from the late 1970's. 
Here, given an abelian group $\Gamma$ and a subset $A\subseteq \Gamma$, we define the \emph{Cayley graph}, $\mathrm{Cay}_\Gamma(A)$ generated by $A$ on $\Gamma$ to be the graph with vertex set $\Gamma$ and containing an arc $\Vec{uv}$ if and only if $v - u \in A$. For a digraph $D$, we denote by $\alpha(D)$ the independence number of its underlying undirected graph obtained by removing all the edge directions.  Our main result is the following Ramsey--Tur\'{a}n variant of Roth's theorem.

\begin{theorem}\label{thm:main}
Let $\calL : \sum_{i\in [k]} c_i x_i = 0$ be a homogeneous linear equation with $k \ge 3$ and $ c_1,\ldots, c_k \in \mathbb{Z} \setminus \{0\}.$ Let further $p$ be a prime. Then, the following are equivalent.
    \begin{enumerate}
        \item[$(a)$] Every solution-free set $A$ of $\mathcal{L}$ in $\mathbb{F}_p$ with $\alpha(\mathrm{Cay}_{\mathbb{F}_p} (A)) = o(p)$ has size $o(p)$. 
        \item[$(b)$] There exists a non-empty set $S\subseteq [k]$ such that $\sum_{s\in S} c_s = 0$.
    \end{enumerate}
\end{theorem}

In particular, our theorem classifies which equations have the property that any solution-free set must either be small or contain a large independent set in its Cayley graph. We note that this graph theoretic property of our set $A$ is natural from the additive point of view as well, it translates to saying that for any $B \subseteq \mathbb{F}_p$ of size $|B| \ge \Omega(p)$ we have $(B-B) \cap A \neq \emptyset.$ This is essentially\footnote{They work with the natural analogue of this definition over $\mathbb{Z}$.} the language used  by Erd\H{o}s and S\'ark\"ozy who refer to this property as $A$ being a difference intersector set.

Obtaining quantitative versions of classical results from additive number theory has attracted a lot of attention in recent years (see e.g.\ \cite{bloom-maynard,kelley-meka,improved-szemeredi}, and the vast number of references therein). Similarly, there are plenty of interesting quantitative questions in classical Ramsey--Tur\'an theory (see e.g.\ \cite{Simonovits-Sos}). With this in mind, as well as to more explicitly state the dependency between the terms hidden behind the $o(p)$ notation in \Cref{thm:main}, we introduce the following definition.

\begin{definition*}
    Let $p$ be a prime number and $\calL$ be a homogeneous linear equation. For a positive real number $\ve > 0$, we denote by $D(\calL, \ve, p)$ the maximum size of a solution-free set $A$ of $\calL$ in $\mathbb{F}_p$ such that $\alpha(\mathrm{Cay}(\mathbb{F}_p, A)) \leq \ve p$. Then, we define
    $$
        d(\calL, \ve) \defeq \limsup_{p\to \infty ~\\~ p: \text{prime}} \frac{D(\calL, \ve, p)}{|\mathbb{F}_p|}.
    $$
\end{definition*}

Armed with this definition, we can rephrase \Cref{thm:main}. Namely, part $(a)$ states that for all $\delta > 0$, there exists $\ve = \ve(\calL, \delta)$ such that $d(\calL, \ve) < \delta$. In fact, if a homogeneous linear equation $\calL$ has $k$ variables (with non-zero coefficients) and satisfies part $(b)$, then our proof of \Cref{thm:main} shows $d(\calL, \ve) < 100^{k+1} k^3\ve$ (see \Cref{thm:degenerate}). We also show that this dependency on $\ve$ can not be improved in general. 

\begin{theorem}\label{thm:schur}
    Let $\ve > 0$ be a real number and $\calL: x + y - z = 0$ be Schur's equation. Then,
    $$
        d(\calL, \ve) = \Theta(\ve).
    $$
\end{theorem}

For an arbitrary equation with a non-zero sum of coefficients, we can show the dependency is polynomial in terms of $\eps$. 

\begin{theorem}\label{thm:poly-lowerbound-intro}
    Let $\calL: c_1x_1 + \ldots +c_k x_k = 0$ be a homogeneous linear equation with $k \ge 3;$ $c_1,\ldots, c_k \in \mathbb{Z} \setminus \{0\}$ such that $c_1+\ldots+c_k \neq 0,$ but there exists $\emptyset \neq S\subseteq [k]: \sum_{s\in S} c_s = 0$. Then, 
    $$
        d(\calL, \ve) = \ve^{\Theta(1)}.
    $$
\end{theorem}

We also raise a number of interesting future directions in \Cref{sec:conc-remarks}.


\section{Preliminaries}\label{sec:prelim}

We first introduce some, mostly standard, notation that we use throughout the paper. For a positive integer $n > 0$, we denote by $[n]$ the set $\{1, \dots, n\}$. As our main arguments rely on graph theory, especially digraphs, we follow the standard notation from (di)graph theory.
Let $D, D'$ be digraphs. We denote by $V(D)$ and $E(D)$ the set of vertices and edges of $D$, respectively. We denote by $D\setminus D'$ the digraph with vertex set $V(D\setminus D') = V(D)$ and edge set $E(D\setminus D') = E(D)\setminus E(D')$. For a vertex $v\in V(D)$, we write $d^+_D(v)$ for the out-degree of $v$, where we frequently omit the subscript if it is clear from context. Similarly, we write $d^-_D(v)$ for the in-degree of $v$. We denote the maximum out and in-degree of $D$, by $\Delta^+(D)$ and $\Delta^-(D)$, respectively. 
Let $U\subseteq V(D)$ be a vertex subset of $V$. Then, we write $D[U]$ for the sub-digraph of $D$ induced on $U$. Finally, let $D_1, \dots, D_k$ be digraphs on a common vertex set $V$. Then we say $(D_1, \dots, D_k)$ is a digraph system.

We will often use the following standard relation between independence number and average degree due to Caro and Wei \cite{caro1979new,wei1981lower}.

\begin{lemma}\label{thm:turan}
    Let $G$ be a graph with average degree $d$. Then $\alpha(G) \geq \frac{|V(G)|}{d + 1}$.
\end{lemma}

\noindent The following lemma is an immediate corollary in our directed setting. 

\begin{lemma}\label{lem:out-degree}
    Let $D$ be a digraph with $\alpha(D) < \frac{|V(D)|}{2r+1}$. Then, $\Delta^+(D) > r$.
\end{lemma}

\begin{proof}
    For the sake of contradiction, assume $\Delta^+(D) \leq r$. This implies that the underlying graph of $D$ contains at most $r|V(D)|$ edges and hence has average degree at most $2r$. So, by \Cref{thm:turan} there is an independent set of size at least $ \frac{|V(D)|}{2r+1}$, contradicting our assumption.
\end{proof}

\noindent A final tweak on \Cref{thm:turan} allows us to remove some edges and retain control of the independence number.

\begin{lemma}\label{lem:remove}
    Let $D, D'$ be $n$-vertex digraphs and let $r \geq 0$ be a non-negative integer. Assume $\Delta^+(D') \leq r$ or $\Delta^-(D') \leq r$. Then, we have
    $$
        \alpha(D\setminus D') \leq (2r + 1) \alpha(D).
    $$
\end{lemma}

\begin{proof}
    Let $I \subseteq V(D)$ be an independent set of $D\setminus D'$. By \Cref{lem:out-degree} applied to $D'[I]$ there exists an independent set $I' \subseteq I$ of $D'$ of size $|I'| \geq \frac{|I|}{2r + 1}$. Since $I'$ is a subset of $I$, it is also an independent set in $D \setminus D'$, so it is an independent set in $D$. Thus, we have 
    $$
        \frac{|I|}{2r + 1} \leq |I'| \leq \alpha(D).
    $$
    From this, we obtain the desired inequality $\alpha(D\setminus D') \leq (2r+1)\alpha(D)$.
\end{proof}


\section{Proof of the classification theorem}\label{sec:proof}
In this section, we prove our classification theorem (\Cref{thm:main}). We refer to a homogeneous linear equation $\calL = \sum_{i\in [k]} c_i x_i$, with $k \ge 3; c_1,\ldots, c_k \in \mathbb{Z} \setminus \{0\}$ as \emph{degenerate} if there exists a non-empty subset $S\subseteq [k]$ such that $\sum_{s\in S} c_s = 0$. 
Then, \Cref{thm:main} can be rephrased as $\lim_{\ve \to 0} d(\calL, \ve) = 0$ if and only if $\calL$ is degenerate.


\subsection{Degenerate equations}\label{subsec:degenerate}

The aim of this section is to prove the implication $(b) \Rightarrow (a)$ of \Cref{thm:main}, which we state more precisely in the following theorem.

\begin{theorem}\label{thm:degenerate}
    Let $\ve > 0$ be a positive real number and $\calL: \sum_{i\in [k]} c_i x_i = 0$ be a degenerate equation. Then, we have
    $$
        d(\calL, \ve) \leq 100^{k+1}k^3\ve.
    $$
\end{theorem}

Before turning to the proof, we show it in a very simple special case of Schur's equation to help build intuition and discuss some immediate generalizations to motivate a number of notions and tools we use in our proof.

A set $X \subseteq \mathbb{F}_p \setminus \{0\}$ not containing a Schur triple is equivalent to  $\mathrm{Cay}_{\mathbb{F}_p}(X)$ not containing a transitive triangle (since $y+z=x \Leftrightarrow z=x-y$). This implies that if $X$ does not contain a Schur triple, then the out-neighborhood of any vertex in $\mathrm{Cay}_{\mathbb{F}_p}(X)$ must be an independent set. On the other hand, all out-neighborhoods have size precisely $|X|$, so $\alpha(\mathrm{Cay}_{\mathbb{F}_p}(X)) \le o(p) \implies |X|\le o(p),$ as desired. There is a small caveat we ignored here so far, namely that we only get to assume that there are no Schur triples with \emph{distinct} entries. This changes the argument only slightly in that within an out-neighborhood every vertex may have at most one out-edge, so we find an independent set of size at least $|X|/3,$ courtesy of \Cref{lem:remove}.

The argument becomes significantly more complex for more involved equations. Firstly, if our equation was of the form $x_1-x_2+x_3+\ldots+x_{k-1}+x_k=0$ very little needs to change, namely we wish to find a rainbow directed path of length $k-2$ (whose colors will play the role of $x_3,\ldots, x_{k}$) with start and end point within an out-neighborhood of some vertex with the edges joining it to them being of different colors from any used on the path (these will be $x_{1}$ and $x_2$). We can find either a linear-sized independent set or such a path within the out-neighborhood itself. Namely, one can iteratively remove vertices with out-degree less than $3k$ (since these vertices have average degree at most $3k$ there cannot be too many of them, or we would find a large independent set by \Cref{thm:turan}). We are left with a graph with minimum out-degree at least $3k$ in which we can greedily build a desired path (when choosing a new vertex we need to avoid already used ones, already used colors and once selected we delete the out-neighbor connected with the color we just used so in total at most $3k-5$ out-neighbors are blocked, so the greedy process can proceed).   

The next level of complexity is well illustrated by equations of the form $x_1-x_2+c_3x_3+\ldots+c_{k-1}x_{k-1}+c_kx_{k}=0.$ Here, we still wish to pick a ``rainbow'' path but the edges should come no longer from $x_3, \dots, x_k$ but from $c_3x_3,\dots, c_k x_k$ while colors are still $x_3,\dots, x_k$. Consequently, at different stages in our path building process, we need to pick an edge from different graphs $\mathrm{Cay}_{\mathbb{F}_p}(c_iX)$ while avoiding a prescribed set of colors, already used vertices and making sure the final and starting vertex connect to our initial vertex using different colors. Since we are not working on the same graph at different stages, some care is needed to ensure this is doable. However, since all these Cayley graphs are isomorphic to $\mathrm{Cay}_{\mathbb{F}_p}(X)$ for which we have control of its independence number, a more careful version of the greedy argument does go through.

The final, major difficulty is that we are not guaranteed two coefficients with opposite signs. Their role is replaced by the subset of the coefficients which sum to zero. Here, we exploit Roth's Theorem (\Cref{thm:roth}) applied repeatedly to obtain a weak supersaturation result to obtain many solutions using mutually disjoint elements of $X$. We then ``give-up'' on the final element in each of these solutions and repeat the above argument inside a suitably translated set of these final elements and find a ``suitable'' path within this set. In order to do this, we first abstract the graph theoretic question lurking in the background, for which the following definitions will come in handy.

Here, we remind the reader that an edge-coloring of a digraph is said to be \emph{proper} if every vertex is incident to at most one out-edge and at most one in-edge of any given color.

\begin{definition*}
    Let $k, \ell > 0$ be positive integers. Let $D_1, \dots, D_k$ be properly edge-colored digraphs with vertex set $V$.
    Let $f$ be a function assigning to each vertex in $V$ a set of up to $\ell$ colors, with $f(v) \cap f(u) = \emptyset$ for any distinct $v,u \in V$. Then, we refer to $(D_1,\dots, D_k; f)$ as an \emph{$\ell$-bounded restricted digraph system}. 
\end{definition*}

In our application, $D_1,\ldots, D_k$ will be Cayley graphs of $\mathbb{F}_p$ generated by $c_i X \subset \mathbb{F}_p,$ so in particular will always be isomorphic. In addition, for each $i$ and any $x \in X$, we give the edge coming from $c_i x$ the color $x$ in $D_i$. The function $f$ encodes certain sets of ``forbidden'' colors. The following definition captures the properties of the path we will need for our argument.

\begin{definition*}
    Let $1 \leq k' \leq k$ be positive integers and $(D_1, \dots, D_k; f)$ be a restricted digraph system on a common vertex set $V$. Then, we say a directed path $P = v_1, \dots, v_{k'+1}$ is a \emph{proper rainbow directed path} if $P$ satisfies the following for any distinct $i, j\in [k']$.
    \begin{enumerate}[label= \upshape\textbf{\Alph{propcounter}\arabic{enumi}}]
        \item $\overrightarrow{v_iv_{i+1}}\in E(D_i)$, \label{def: rainbow 1}
        \item the color of $\overrightarrow{v_iv_{i+1}}$ in $D_i$ is distinct from the color of $\overrightarrow{v_jv_{j+1}}$ in $D_j$, and \label{def: rainbow 2}
        \item the color of $\overrightarrow{v_iv_{i+1}}$ in $D_i$ is not in $f(v_1)\cup f(v_{k'+1})$. \label{def: rainbow 3}
    \end{enumerate}\stepcounter{propcounter}
 \end{definition*}

The following lemma is the heart of our proof of \Cref{thm:degenerate}. It guarantees the existence of a proper rainbow directed path in an $\ell$-bounded restricted digraph system for which we have some control of the independence number of each digraph. 

\begin{lemma}\label{lem:rainbow-path}
    Let $1 \leq k' \leq k$ be positive integers.
    Let $(D_1, \dots, D_k; f)$ be an $\ell$-bounded restricted digraph system on a common vertex set $V$. If $\alpha(D_i) \leq \frac{|V|}{100^{k'}\ell^2}$ for all $i\in [k'],$ then $(D_1, \dots, D_k; f)$ contains a proper rainbow directed path of length $k'$.
\end{lemma}

\begin{proof}
    We use induction on $k'$. For the base case, suppose $k' = 1$ and $\alpha(D_1) \leq \frac{|V|}{100\ell^2}$. We let $D'$ be a subdigraph of $D_1$ consisting of edges $\overrightarrow{vu}$ with color in $D_1$ belonging to $f(u)$. 
    Since the colorings are assumed to be proper, we have $\Delta^-(D') \leq \max_{u\in V}|f(u)| \leq \ell$. Thus, by \Cref{lem:remove},
    \begin{equation}\label{eq:D'}
        \alpha(D_1 \setminus D') \leq (2\ell+1) \alpha(D_1) \leq \frac{|V|}{30\ell}.
    \end{equation}
    This combined with \Cref{lem:out-degree} implies that there exists a vertex $v\in V$ such that $d^+_{D_1\setminus D'}(v) \geq \ell+1$. This allows us to choose an edge $\overrightarrow{vu} \in E(D_1\setminus D')$ with color in $D_1$ not belonging to $f(v)$. This edge is a desired proper rainbow directed path. Indeed, properties a) and b) hold trivially, and for c) we ensured the color of $\overrightarrow{vu}$ is not in $f(v)$ when picking $u$ and that it is not in $f(u)$ by ensuring $\overrightarrow{vu} \notin E(D')$. This proves the base case, $k' = 1$.

    Now assume that $k' = r+1$ and that the statement of \Cref{lem:rainbow-path} holds for all $k' \leq r$.
    Suppose now that for each $i\in [r+1]$, we have $\alpha(D_i) \leq \frac{|V|}{100^{r+1} \ell^2} $. Let $F$ be an auxiliary digraph on $V$ such that $\overrightarrow{vu}\in E(F)$ if and only if there is a proper rainbow directed path of length $r$ from $v$ to $u$. Then, by applying the induction hypothesis to the subsystem induced on any $|V|/100$ vertices of $V$ we conclude $F$ contains an edge within any vertex set of size at least $|V|/100$. In other words,
    \begin{equation}\label{eq:F}
        \alpha(F)\leq \frac{|V|}{100}.
    \end{equation}
    
    \noindent Let $U\defeq \{v\in V: d_F^-(v) \geq 1\}$ and note that $V\setminus U$ is an independent set of $F$, so 
    \begin{equation}\label{eq:U}
        |U| \geq \frac{99|V|}{100}.
    \end{equation}

    \noindent Let $D'$, similarly as before, be a subdigraph of $D_{r+1}$ 
    consisting of edges $\overrightarrow{vu}$ with color in $D_{r+1}$ belonging to $f(u)$. 
    Again, by the assumption that the coloring of $D_{r+1}$ is proper, we have $\Delta^-(D') \leq \max_{u\in V}|f(u)| \leq \ell$. Thus, by \Cref{lem:remove},
    \begin{equation}\label{eq:D'r}
        \alpha(D_{r+1} \setminus D') \leq (2\ell+1) \alpha(D_{r+1}) \leq \frac{3|V|}{100^{r+1}\ell}.
    \end{equation}

    \noindent Let $F'$ be the subdigraph of $D_{r+1}\setminus D'$ induced on $U$. Then, by \eqref{eq:D'r} and \eqref{eq:U}, we have 
    \begin{equation}\label{eq:F'}
        \alpha(F') \leq \alpha(D_{r+1} \setminus D) \leq \frac{3|V|}{100^{r+1}\ell} \leq \frac{|U|}{3\cdot 100^r \ell}.
    \end{equation}

    By \Cref{lem:out-degree} and \eqref{eq:F'}, we conclude that there is a vertex $v_{r+1}\in U$ such that $d^+_{F'}(v_{r+1}) \geq 100^r \ell > \ell +  3r + 1$. Since $v_{r+1}\in U$ there exists $v_1$ such that $ \overrightarrow{v_1v_{r+1}} \in F$ so in particular there is a proper rainbow directed path $v_1,v_2,\ldots, v_{r+1}$ (so of length $r$). Since $d^+_{F'}(v_{r+1}) > \ell + 3r + 1$ we can find an edge $\overrightarrow{v_{r+1}v_{r+2}} \in E(D_{r+1}\setminus D')$ such that all the following properties hold. 
    \begin{itemize}
        \item $v_{r+2}$ is distinct from $v_1,\ldots, v_r,v_{r+1}$.
        \item The color of $\overrightarrow{v_{r+1}v_{r+2}}$ in $D_{r+1}$ is distinct from the color of $\overrightarrow{v_{i}v_{i+1}}$ in $D_{i+1}$ for any $1\le  i \le r$.
        \item The color of $\overrightarrow{v_{r+1}v_{r+2}}$ in $D_{r+1}$ is not in $f(v_1)$.
        \item The color of $\overrightarrow{v_{i}v_{i+1}}$ in $D_{i+1}$ is not in $f(v_{r+2}),$ for any $1 \le i \le r$. 
    \end{itemize}
    The first bullet excludes up to $r+1$, the second up to $r,$ the third up to $\ell$, and the last one up to $r$ choices for $v_{r+2}$, so we can indeed make such a choice. Here, for the last point, we used that $f$ assigns disjoint sets of colors to distinct vertices, which guarantees that any color already appearing on the path blocks at most one choice for $v_{r+2}$.
    
    We claim that $v_1,\ldots, v_{r+1},v_{r+2}$ is a desired proper rainbow directed path of length $r+1$. Indeed, using that $v_1,v_2,\ldots, v_{r+1}$ is a proper rainbow path combined with $\overrightarrow{v_{r+1}v_{r+2}} \in E(D_{r+1})$ ensures property a), combined with the first bullet ensures it is a directed path, combined with the second bullet ensures property b), and combined with the third bullet above it ensures colors in $f(v_1)$ are avoided. Finally, since $\overrightarrow{v_{r+1}v_{r+2}} \notin E(D')$ we know that the color of $\overrightarrow{v_{r+1}v_{r+2}}$ in $D_{r+1}$ is not in $f(v_{r+2})$ and the fourth bullet above ensures this holds for all the colors of $\overrightarrow{v_{i}v_{i+1}}$ in $D_{i+1}$ for $i\le r$ as well, establishing property c).
    This completes the induction step and with it the proof.
\end{proof}

\noindent We are now ready to prove \Cref{thm:degenerate}.

\begin{proof}[Proof of \Cref{thm:degenerate}]
    Let $p$ be a sufficiently large prime number so that \Cref{thm:roth} guarantees that any subset of $\mathbb{F}_p$ of size at least $\eps p$ contains a solution of any homogenous linear equation with $3\le k' \le k$ variables with coefficients among $\{c_1,\ldots, c_k\}$ summing to zero.
    
    Let $A$ be a subset of $\mathbb{F}_p$ that satisfies the following.
    \begin{enumerate}[label= \upshape\textbf{\Alph{propcounter}\arabic{enumi}}]
        \item $|A| > 100^{k+1} k^3 \cdot \ve p$, and \label{B1}
        \item $\alpha(\mathrm{Cay}_{\mathbb{F}_p}(A)) \leq \ve p$. \label{B2}
    \end{enumerate}\stepcounter{propcounter}
    Then, it suffices to show that $A$ contains a solution of $\calL$ using distinct elements of $A$.
    
    For each $i\in [k]$, let $D_i\defeq \mathrm{Cay}_{\mathbb{F}_p}(c_i\cdot A)$. As $c_i$ is a fixed non-zero integer, and $p$ is a sufficiently large prime number, we note that $D_i$ is isomorphic to $\mathrm{Cay}_{\mathbb{F}_p}(A)$. Thus, for each $i\in [k]$ and $v\in \mathbb{F}_p$, we have $d^+_{D_i}(v) = d^-_{D_i} (v)= |A| > 100^{k+1} k^3 \cdot \ve p$ and $\alpha(D_i) \leq \ve p$.
    We also assign a proper edge coloring to each $D_i$, where an edge of $D_i$ arising from $c_i a$ with $a \in A$ gets color $a$ (note here that for $c_i \neq 1$ this is \emph{not} the usual Cayley coloring, but rather a ``translation'' of it). 
    
    Since $\calL$ is degenerate, there is a non-empty set of indices $S\subseteq [k]$ such that $\sum_{s\in [S]} c_s = 0$. Without loss of generality, let $S = [k']$ for some integer $2 \leq k' \leq k$. If $k' = k$, then by \Cref{thm:roth}, we can find a desired solution and are done. Thus, we may assume $k' < k$. 
    
    Let $\calL':= \sum_{i\in [k']} c_i x_i = 0$. Let us for the moment, assume $k'\ge 3$. Then, by \Cref{thm:roth}, every subset of $A$ of size at least $\eps p$ in $\mathbb{F}_p$ contains a solution to $\calL'$ using distinct elements. By deleting these elements from $A$ and repeating until we are left with less than $\eps p$ elements in $A$, we can find a set of solutions $\calS = \{\mathbf{x}_1, \dots, \mathbf{x}_{\ell}\} \subseteq A^{k'}$ of $\calL'$ such that
    \begin{enumerate}[label= \upshape\textbf{\Alph{propcounter}\arabic{enumi}}]
        \item $\ell \geq 100^{k} k^2  \cdot \ve p$,
        \item for $(r, i) \neq (s, j) \in [\ell]\times [k']$,
     if $r\neq s$ or $k'\notin \{i,j\}$, then $(\mathbf{x}_r)_i \neq (\mathbf{x}_s)_j$. 
    \end{enumerate}\stepcounter{propcounter}
    If $k'=2$, we must have $c_1=-c_2$, so any $x_1=x_2$ makes a solution (albeit not one using distinct elements\footnote{Which is the main reason for the insistence on $k \ge 3$ throughout the paper.}). Here, we can also trivially find a set of solutions as above.\footnote{The condition ``$r\neq s$ or $k'\notin \{i,j\}$" on (b) is to ensure that it is satisfied for $k'=2$.}

    Denote by $U$ the set of $(-c_{k'})$-dilations of $k'$-th coordinates of tuples in $\calS$, that is $$U\defeq  \{-c_{k'}(\mathbf{x}_1)_{k'},-c_{k'}(\mathbf{x}_2)_{k'},\ldots, -c_{k'}(\mathbf{x}_\ell)_{k'}\} \subseteq \mathbb{F}_p.$$ Then, we have
    \begin{equation}\label{eq:V-size}
        |U| = \ell \geq 100^{k} k^2 \cdot \ve p.
    \end{equation}

    Next, we define $f: U\to 2^{A}$ by setting $f(u) := \{(\mathbf{x}_r)_1, \dots, (\mathbf{x}_r)_{k'}\},$ where $r$ is the index of the solution for which $u$ is the ($-c_{k'}$)-dilation of its last coordinate (i.e.\ $u = - c_{k'}(\mathbf{x}_r)_{k'} \in U$). We note that such $r$ is unique since the last coordinates are all mutually different (and multiplied by the same non-zero factor $c_{k'}$), so $f$ indeed assigns disjoint sets of colors to distinct vertices in $U$. So, $(D_{k'+1}[U], \ldots, D_k[U]; f)$ is a $k'$-bounded restricted digraph system. Furthermore, by \eqref{eq:V-size} and \ref{B2}, we have that for all $k'+1 \leq i \leq k$, 
    $$
        \alpha(D_i[U])\le \alpha(D_i)=\alpha(\mathrm{Cay}_{\mathbb{F}_p}( A)) \leq \frac{|U|}{100^k k^2}.
    $$
    Thus, by \Cref{lem:rainbow-path}, the restricted digraph system $(D_{k'+1}[U], \dots, D_k[U]; f)$ contains a proper rainbow directed path $P = v_{k'},v_{k'+1},\dots,v_{k}$ of length $k - k'$.
    Let $v_{k'} = -c_{k'} (\mathbf{x}_r)_{k'}$ for some $r\in [\ell]$ (which exists since $v_{k'} \in U$). For each $i\in [k'-1]$, let $z_i \defeq (\mathbf{x}_r)_i$. Then,
    \begin{equation}\label{eq:first}
        \sum_{i\in [k'-1]} c_i z_i = v_{k'},
    \end{equation}
    and all $z_i$ for $i \le k'-1$ are distinct by property (b)\footnote{Note that we are \emph{not} using the last coordinate here.}.

    Now, for all $k'+1 \leq j \leq k$, let $z_j \defeq c_j^{-1}\cdot (v_{j}-v_{j-1})$. We note that, by property \ref{def: rainbow 1} of a proper rainbow directed path, we have $v_{j} - v_{j-1} \in E(D_j)$ which implies $z_j \in A,$ for all $j \ge k'+1.$ By property \ref{def: rainbow 2}, we also know that $z_j$ are all distinct for $j \ge k'+1$. Since $P$ is a directed path, we also have
    \begin{equation}\label{eq:second}
        \sum_{k'+1 \leq j \leq k} c_j z_j = v_k - v_{k'}. 
    \end{equation}

    Finally, let $z_{k'} \defeq c_{k'}^{-1}\cdot(-v_k)$, so that
    \begin{equation}\label{eq:third}
        c_{k'}z_{k'} = -v_k.
    \end{equation}
    Note that since $v_k \in U,$ by the definition of $U$, the element $z_{k'}\in A$.
    Then, by combining \eqref{eq:first}--\eqref{eq:third}, $z_1,\ldots,z_k$ make a solution for $\calL$ in $A$, i.e.\

    $$
        \sum_{i\in [k]} c_i z_i = 0.
    $$
    Finally, to ensure all the $z_i$ are distinct note that $f(v_{k'})=\{(\mathbf{x}_r)_{1},\ldots, (\mathbf{x}_r)_{k'}\}\supseteq \{z_1,\ldots, z_{k'-1}\}$ and $z_{k'} \in f(v_k)$. So, indeed, since $z_{k'+1},\ldots, z_k$ are the colors used by the edges of our proper rainbow directed path, we conclude by its properties \ref{def: rainbow 2} and \ref{def: rainbow 3} that all $z_i$ are distinct.
 This completes the proof.    
\end{proof}


\subsection{Non-degenerate equations}\label{subsec:non-degenerate}

In this section, we prove the other direction $(a) \Rightarrow (b)$ of \Cref{thm:main}. The following theorem directly implies this.

\begin{theorem}\label{thm:non-degenerate}
For a non-degenerate homogeneous linear equation $\calL = \sum_{i\in [k]} c_i x_i = 0$,  there exists $\beta > 0$ and $p_0 > 0$ such that the following holds. For all prime numbers $p > p_0$, there exists a solution-free subset $A\subseteq \mathbb{F}_p$ of $\calL$ that satisfies the following.
\begin{enumerate}
    \item[$\bullet$] $|A| \geq \beta p,$ and
    \item[$\bullet$] $\alpha(\mathrm{Cay}_{\mathbb{F}_p}(A)) = O_{\calL} (\frac{p}{\log p})$. 
\end{enumerate}
In particular, we have $d(\calL, \ve) \geq \beta > 0$ for all $\ve > 0$.
\end{theorem}

\begin{proof}[Proof of \Cref{thm:non-degenerate}]
    Choose $t \in \ZZ$ such that $t > \sum_{i=1}^k |c_i| := C$. Define $Z = X \cup Y \subseteq \ZZ$ where
    $$
        X = \left\{rt + 1: 0 \leq r \leq \frac{p}{2Ct}\right\},
    $$
    $$
        Y = \{t^z - t^w : 1 \leq w < z \leq \log_t \sqrt{p}\}
    $$
    For each $a\in Z$, we denote by $\bar{a}$ the element in $\mathbb{F}_p$ such that $\bar{a} \equiv a \pmod{p}$.
    Let $A\defeq \{\bar{a}\in \mathbb{F}_p: a\in Z\}$ and let $\beta = \frac{1}{2Ct}$.
    We will show that $A$ is the desired subset.
    Note that any distinct $z_1, z_2 \in Z$, satisfy $0 < |z_1 - z_2| < p$, so $|A| = |Z| \geq |X| \geq \beta p$ holds. 
    
    \begin{claim*}\label{clm:solution-free}
        $A$ is a solution-free subset of $\calL$ in $\mathbb{F}_p$. 
    \end{claim*}
    
    \begin{claimproof}
        Suppose towards a contradiction that $A$ contains a solution $(\bar{a}_1,\ldots, \bar{a}_k)$ such that $\sum_{i=1}^k c_i \bar{a}_i = 0$ over $\mathbb{F}_p,$ and $\bar{a}_i$ are all distinct. Note that for large enough $p$, $\max(X) < \frac{p}{C},$ and $\max(Y) \leq t^{\log_t \sqrt{p}} = \sqrt{p} < \frac{p}{C}$. Hence, when viewed as an equation over the integers,
        $$
            0\leq \left|\sum_{i=1}^k c_i a_i\right| \leq \sum_{i=1}^k |c_i| |a_i| < \frac{p}{C}\sum_{i=1}^k |c_i| = p.
        $$
        Therefore, we must have that $\sum_{i=1}^k c_i a_i = 0$ over $\ZZ$ as well. 
        
        We now consider two cases: whether or not some of the $a_i$ are contained in $X$. For the first case, assume at least one of $a_i$ lives in $X$. Without loss of generality, for some $1 \leq k' \leq k$, assume that $a_i \in X$ for $i \leq k'$ and $a_j \in Y$ for $k' < j \leq k$.
        Note that $\sum_{i > k'} c_i a_i$ is divisible by $t$ since every $a_i \in Y$ is of the form $t^z - t^w$ which is divisible by $t$. This implies that $\sum_{i=1}^{k'} c_i a_i$ is divisible by $t$ as well.
        Since every $a_i \in X$ satisfies $a_i \equiv 1 \pmod{t}$, it follows that $\sum_{i=1}^{k'} c_i \equiv 0 \pmod{t}$. By our choice of $t$, we have $\sum_{i=1}^{k'} |c_i| < t$. These observations imply the number $\sum_{i=1}^{k'} c_i$ is equal to $0$, contradicting the fact that $\calL$ is non-degenerate.

        Thus, we may assume that all of $a_i$ are contained in $Y$. Write $a_i = t^{z_i} - t^{w_i}$ for each $i \in [k]$. Without loss of generality, assume $y_1$ is the largest among the $y_i$'s. Then, we have $z_1 \geq z_i > w_i$ for all $i \in [k]$. Let $I \defeq \{i \in [k]: z_i = z_1\}$. Then, we have 
        $$
            \sum_{i \in I} c_i t^{z_1} = \sum_{i \in I} c_i t^{w_i} - \sum_{j \in [k]\setminus I} c_j(t^{z_j}-t^{w_j})
        $$
        Since $\sum_{i \in I} c_i \neq 0$ by assumption, the inequality
        \begin{equation}\label{eq:t}
            \left|\sum_{i \in I} c_i t^{z_1} \right| \geq t^{z_1}.
        \end{equation} 
        holds.
        On the other hand, by our choice of $t$, 
        $$
            \left|\sum_{i \in I} c_i t^{w_i} - \sum_{j \in [k]\setminus I} c_j(t^{z_j}-t^{w_j})\right| \leq \sum_{i=1}^k |c_i| t^{z_1 - 1} < t^{z_1}.
        $$ This contradicts \eqref{eq:t}. Thus, $A$ is indeed a solution-free subset of $\calL$ in $\mathbb{F}_p$. 
    \end{claimproof}

    We now claim that $\mathrm{Cay}_{\mathbb{F}_p}(A)$ has a small independence number.

    \begin{claim*}
        $\alpha(\mathrm{Cay}_{\mathbb{F}_p}(A)) = O_t(\frac{p}{\log p})$.
    \end{claim*}
    
    \begin{claimproof}
        Let $G$ be the underlying graph of $\mathrm{Cay}_{\mathbb{F}_p}(A)$. We note that $A$ is a vertex-transitive graph, hence the inequality
        \begin{equation}\label{eq:clique-number}
            \alpha(G) = \frac{p}{\chi_f(G)} \leq \frac{p}{\omega(G)}
        \end{equation} 
        holds, where $\chi_f(G)$ and $\omega(G)$ are the fractional chromatic number and the clique number of $G$, respectively.
        Consider the set $B \defeq \{\bar{t^i}: 1 \leq i \leq \log_t \sqrt{p}\}$. Then, for each $1\leq i < j \leq \log_t \sqrt{p}$, we have $t^j - t^i \in A$ because of the definition of $X$. Hence, $(\bar{t^i}, \bar{t^j})$ is an edge in $G$ for all $1 \leq i < j \leq \log_t(\sqrt{p})$. Therefore, $B$ forms a clique of size at least $\Omega_t(\log p)$ in $G$. By \eqref{eq:clique-number}, we have $\alpha(G) = O_t(\frac{p}{\log p})$.
    \end{claimproof}

    \noindent From the above two claims, it follows that $A$ indeed satisfies the desired properties.
\end{proof}


\section{Quantitative lower bounds}\label{sec:proof-schur-lowerbound}

In this section, we prove our quantitative lower bounds on $d(\calL, \ve)$. We start with the lower bound for \Cref{thm:schur}. As an ingredient in the proof, we will use the lower bound on Ramsey numbers of a triangle vs large independent sets.
The correct asymptotic upper bound here was determined by Ajtai, Koml\'os and Szemer\'edi \cite{AKS} with the best known constant due to Shearer \cite{shearer}. The correct asymptotic lower bound was first proved by Kim \cite{kim-triangle-ramsey} with a number of remarkable papers (e.g.\ \cite{triangle-free-1,triangle-free-2}) improving the constant factor leading to the current state of the art bound from a very recent paper \cite{new-ramsey} bringing the bound to within a $3+o(1)$ factor of Shearer's bound. Any correct asymptotic lower bound would suffice for our purposes, but we pick the following for clarity of presentation.  

\begin{theorem}\label{thm:ramsey}
    For any sufficiently large integer $r$ there exists a triangle-free graph $G$ with $\alpha(G)=r$ and at least $\frac{r^2}{8\log r}$ vertices.
\end{theorem}

We are now ready to prove our lower bound for Schur's equation. At a high level we will embed a Ramsey graph provided by the above theorem onto a very small portion of our Cayley graph in order to ensure the independence number is not too large, we will then add to it an intersection of a ``middle'' interval with said independent set (which we want to ensure is not too small either) to ensure our set is reasonably large while still avoiding any Schur triples.  

\begin{theorem}\label{thm:schur-lowerbound}
    Let $1\ge \ve > 0$ be a real number and $\calL: x + y - z = 0$ be Schur's equation. Then,
    $$
        d(\calL, \ve) = \Omega(\ve).
    $$
\end{theorem}
\begin{proof}
Note here that we may assume $\eps<\eps_0$ for any fixed constant $\eps_0>0$ since  $d(\calL, \ve)$ is increasing in $\ve$. Let $p$ be a prime sufficiently larger than $1/\eps$, taking $p>20\cdot (100/\eps)^{2500/\eps^2}$ would suffice.

Pick $t=100/\eps$ and let $G$ be an auxiliary triangle-free graph with $\alpha(G)= t \log t$ and $q=\frac18 t^2 \log t$ vertices guaranteed by \Cref{thm:ramsey} (where we delete superfluous vertices if necessary). Let us also identify the vertex set of $G$ with $[q]$.
Let $X:= \{ 4^j - 4^{\ell}: j\ell \in E(G), j>\ell \}$ and observe that if we set $Q:=\{4,4^2,\ldots, 4^{q}\}$ the subgraph of $\mathrm{Cay}_{\mathbb{F}_p}(X)$ induced on $Q$ is isomorphic to $G$ (using that $p>2 \cdot 4^{q}$ to avoid ``overflows'').
Moreover, this remains true when considering the subgraph induced on $x+Q$ for any $x \in \mathbb{F}_p$.

Next, we claim that $\alpha(\mathrm{Cay}_{\mathbb{F}_p}(X)) \le 10p/t.$ Indeed, if this graph contains a larger independent set $I$, then 
$$ \sum_{x \in \mathbb{F}_p} |(x+Q)\cap I|=|Q||I| > q \cdot 10p/t.$$ This implies that there exists $x \in \mathbb{F}_p$ with $|(x+Q)\cap I| > 10q/t >t \log t$, a contradiction since $(x+Q)\cap I$ is an independent set of $\mathrm{Cay}_{\mathbb{F}_p}(X)$ induced on $x+Q$.

Let us enumerate the elements of $X$ as $x_1,\dots, x_{e}$, where $e=|E(G)|$, and set $X_i:= \{x_1,\ldots, x_i\}$.
Consider the sequence of Cayley graphs $\mathrm{Cay}_{\mathbb{F}_p}(X_0),\ldots, \mathrm{Cay}_{\mathbb{F}_p}(X_e)$. Let $H_i$ be the subgraph of $\mathrm{Cay}_{\mathbb{F}_p}(X_i)$ induced on the interval $[p/3,4p/9]$. Observe that $\alpha(H_0) \ge p/9,$ that $\alpha(H_e) \le 10p/t,$ and that $\alpha (H_i) \ge \alpha(H_{i+1})\ge \alpha(H_i)/2$ since $H_{i+1}$ differs from $H_i$ by an addition of a path forest (consisting of edges with color $x_{i+1}$ induced on $[p/3,4p/9]$) any independent set of $H_i$ contains an independent set of $H_{i+1}$ of at least a half of its size. Let us hence pick an $i$ for which $10p/t\ge \alpha(H_i) \ge 5p/t$. Note that $\alpha(H_i) \le 10p/t$ implies the same bound holds for any subgraph of $\mathrm{Cay}_{\mathbb{F}_p}(X_i)$ induced on an interval of width $p/9$ so in particular $\alpha(\mathrm{Cay}_{\mathbb{F}_p}(X_i))\le 100p/t.$ 

Let $Y$ be an independent set of $H_i$ with $|Y|\geq 5p/t$. 
We take $X_i\cup Y$ as our final generating set, then it has size at least $5p/t \geq \ve p/20$. Moreover the independence number of $\mathrm{Cay}_{\mathbb{F}_p}(X_i\cup Y)$ is smaller than the independence number of $\mathrm{Cay}_{\mathbb{F}_p}(X_i)$, which is at most $100p/t \leq \ve n$. 

Finally, it remains to verify that  $X_i\cup Y$ has no Schur triples. Since $Y \subseteq [p/3,4p/9]$, there are no Schur triples with all elements in $Y$. Furthermore, $Y+Y \subseteq [2p/3,8p/9]$ and since $\max X_i <p/9$ we have $(Y+Y) \cap X_i = \emptyset$. Combined with $(Y-Y) \cap X_i=\emptyset$, which holds since $Y$ is an independent set in $\mathrm{Cay}_{\mathbb{F}_p}(X_i)$, we conclude there are no Schur triples with precisely two elements from $Y$. Since $X_i \pm X_i \subseteq (8p/9,p) \cup [0,p/9)$ which is disjoint from $Y$ we conclude having precisely one element from $Y$ is impossible. Finally, since $X_i \subseteq X$ any Schur triple in $X_i$ would take form $4^{b}-4^{a}=4^{b'}-4^{a'}+4^{b''}-4^{a''}$ with $a>b, a'>b',a''>b''$ and $ba,b'a',b''a'' \in E(G)$. Viewing this equality in base four implies $\{b,a',a''\}=\{a,b',b''\}$ which in turn implies $a'=b''$ or $b'=a''$. By symmetry, we may w.l.o.g.\ assume $a'=b''$, which in turn gives $b=b'$ and $a=a''$. However, this gives rise to a triangle $ba,ba',a'a$ in $G$ which does not exist.
\end{proof}

For our general lower bound, the following result is going to play the role of \Cref{thm:ramsey} in the case of Schur's equation. It is due to Osthus and Taraz \cite{osthus}\footnote{We note that while they state their Theorem 5 in terms of the chromatic number, in the proof they establish the slightly stronger version of the result we quote here in terms of the independence number.} and has a simple proof based on the alterations method. We note that for $\ell=3$ the polylogarithmic factor is a bit weaker when compared to \Cref{thm:ramsey} and using the subsequent famous developments on the $H$-free process one can likely improve this factor here, although it remains open whether even the polynomial term in this result is tight for $\ell \ge 4$. 

\begin{theorem}\label{thm:ramsey-girth}
    For any integer $\ell \ge 3$, there exists $C=C(\ell)>0$ such that for any $n$ there exists an $n$-vertex graph $G$ with no cycles of length at most $\ell$ and $\alpha(G) \le C n^{1-1/(\ell-1)} \log n$.
\end{theorem}

We are now ready to prove our general lower bound. Our strategy is similar as in \Cref{thm:schur-lowerbound} with some key changes.

\begin{theorem}\label{thm:poly-lowerbound}
    Let $\calL: c_1x_1 + \ldots +c_k x_k = 0$ be a homogeneous linear equation with $k \ge 3;$ $c_1,\ldots, c_k \in \mathbb{Z} \setminus \{0\}$ such that $c_1+\ldots+c_k \neq 0$. There exists $C=C(\calL)>0$ such that for any $1> \ve > 0$ we have
    $$
        d(\calL, \ve) \ge \ve^{C}.
    $$
\end{theorem}

\begin{proof}  
Let $M=\max |c_i|$ and pick a prime $r$ between $k M$ and $2kM$. We note that since our equation is fixed, one should treat $r$ as a large constant throughout the proof.
By adjusting $C$ if necessary, we may assume $\eps<\eps_0$ for any fixed constant $\eps_0>0$ since  $d(\calL, \ve)$ is increasing in $\ve$. Let $p$ be a prime sufficiently larger than $1/\eps$.

Let $t=1/\eps$ and let $G$ be an auxiliary graph with no cycles of length at most $k+1,$ with $\alpha(G)\leq  t^{k-1}$ and $q=t^{k}$ vertices guaranteed by \Cref{thm:ramsey-girth} (where we are using the assumption that $t$ is large enough and are being a bit wasteful to simplify numbers). Let us identify the vertex set of $G$ with $[q]$.
Let $X:= \{ r^j - r^{\ell}: j\ell \in E(G), j>\ell \}$ and observe that if we set $Q:=\{r,r^2,\ldots, r^{q}\}$ the subgraph of $\mathrm{Cay}_{\mathbb{F}_p}(X)$ induced on $Q$ is isomorphic to $G$.
Moreover, this remains true when considering the subgraph induced on $x+Q$ for any $x \in \mathbb{F}_p$.

Next, we claim that $\alpha(\mathrm{Cay}_{\mathbb{F}_p}(X)) \le p/t.$ Indeed, if this graph contains a larger independent set $I$, then $$ \sum_{x \in \mathbb{F}_p} |(x+Q)\cap I|=|Q||I| > q \cdot p/t.$$ This implies that there exists $x \in \mathbb{F}_p$ with $|(x+Q)\cap I| > q/t \geq t^{k-1}$, a contradiction since $(x+Q)\cap I$ is an independent set of $\mathrm{Cay}_{\mathbb{F}_p}(X)$ induced on $x+Q$.

Next, let us pick $r'$ to be the smallest prime not dividing any non-zero element in $X_{\Sigma}:=\{\pm x_1 \pm \ldots \pm x_m \mid x_i \in X, m \le r\}\setminus\{0\}$. Note that $X_{\Sigma}$ consists of at most $rq^r$ elements, each of which has absolute value at most $r^{q+1}$. Hence, the absolute value of the product of these elements is at most $r^{2rq^{r+1}}$. On the other hand, the product of the first $2r q^{r+1} \log r$ primes is strictly larger than this implying that $r'$ is among these primes and in particular by the prime number theorem that $r'\le 200 q^{r+2} \cdot r^2 \log^2 r \le t^{2rk} \leq t^{4k^2M}.$
We let $Y$ consist of the elements in $\left [\frac p{2r}- \frac {p}{4r^2},\frac {p}{2r} + \frac {p}{4r^2}\right]$ which are divisible by $r'$. 

We take $X\cup Y$ as our final generating set. It has size at least $|Y|\ge \frac{p}{4r^2r'} \geq \ve^{4 k^2 M+1} p $, where we use that we may assume $1/\ve\geq \frac{1}{4r^2}$. 
Moreover, $\alpha(\mathrm{Cay}_{\mathbb{F}_p}(X\cup Y)) \le \alpha(\mathrm{Cay}_{\mathbb{F}_p}(X))\le p/t=\eps p.$ 

Finally, it remains to verify that  $X\cup Y$ is solution-free for $\calL$. 
Suppose, towards a contradiction that we have a solution $(x_1,\ldots,x_k)$ in $(X \cup Y)^k$ which after relabeling coefficients of $\calL$ w.l.o.g.\ satisfies $c_1x_1+\ldots+c_k x_k=0$ with $x_1,\ldots,x_i \in X$ and $x_{i+1},\ldots, x_k \in Y$ for some $0\le i \le k$. Note that $|c_1x_1+\ldots+c_i x_i| < r^{q+1} \le p/(4r)$. Notice next that $c_{i+1}x_{i+1}+\ldots+c_kx_k \in \left[(c_{i+1}+\ldots+c_k)\cdot \frac p{2r}- r\cdot \frac {p}{4r^2},(c_{i+1}+\ldots+c_k)\cdot \frac p{2r}+ r\cdot \frac {p}{4r^2}\right],$ since $|c_{i+1}|+\ldots+|c_k| \le r$. These two observations imply that $c_{i+1}+\ldots+c_k=0$ and hence that $c_1+\ldots +c_i \neq 0$, so in particular $i \ge 1$. Since $Y$ consists only of elements divisible by $r'$ we know that $r' \mid c_1x_1+\ldots+c_i x_i$, and since this is in $X_{\Sigma}$, by our choice of $r'$, we conclude $c_1x_1+\ldots+c_i x_i=0$. 
Since $x_1,\ldots, x_i \in X$ this implies there exist $1\le a_j <b_j \le q$ with $x_j=r^{b_j}-r^{a_j},$ where $\{a_j,b_j\} \in E(G)$ are distinct\footnote{We may also ensure this by (partially) canceling out occurrences of the same choice for $x_j$ with $c_i$'s of opposite sign.} for all $j$. This gives $c_1r^{b_1}+\ldots +c_ir^{b_i}=c_1r^{a_1}+\ldots+c_ir^{a_i}.$ By viewing this equation in base $r$, since $|c_1|+\ldots+|c_i|<r$ we can conclude that any $x \in \{a_1,b_1,\ldots, a_i,b_i\}$ appears in at least two distinct $\{a_j,b_j\}$'s. This means that the subgraph of $G$ induced on $\{a_1,b_1,\ldots, a_i,b_i\}$ has minimum degree at least two and hence contains a cycle of length at most $i \le k$, which is a contradiction to how we picked our $G$.
\end{proof}


\section{Concluding remarks}\label{sec:conc-remarks}

In this paper, we solve a natural Ramsey--Tur\'an analogue of Roth's classical theorem on density regularity of equations in $\mathbb{F}_p$. There are many interesting future directions to explore, some of which we discuss briefly below. 

\noindent \textbf{Other groups. } Perhaps the most immediate further question is: what happens for groups other than $\mathbb{F}_p$? By adapting our definitions and proofs in an appropriate way, one can show an appropriate analogue of  Theorem~\ref{thm:main} over $\mathbb{Z}$. There are a few minor caveats, but such translations seem fairly standard in additive combinatorics, so we leave this as a potentially nice undergraduate student project. For the general case of finite abelian (or even non-abelian) groups $\Gamma$ and equation $\mathcal{L}: \sum_{i\in [k]} c_i x_i=0$, there are certain difficulties in adapting our proofs. For degenerate equations, one issue is that some $c_i$ may not be coprime with $|\Gamma|$, so that $c_i\Gamma$ may not be isomorphic to $\Gamma$. For non-degenerate equations, it even seems difficult to find a construction when there exists $I\subseteq [k]$ such that $\sum_{i\in I} c_i$ is not coprime with $|\Gamma|$.

\noindent \textbf{Quantitative dependence. } Our proof shows that the value of $d(\mathcal{L},\ve)$  is at least $\ve^{O(k^2 M)}$, where $k$ is the number of variables and $M$ is the maximum absolute value of the coefficients in $\mathcal{L}$. We determine the right exponent of $\ve$ in the case of Schur's equation, and it would be interesting to obtain such a result in general. We conjecture the following strengthening of our \Cref{thm:poly-lowerbound-intro}.
\begin{conjecture}
   Let $\calL: c_1x_1 + \ldots +c_k x_k = 0$ be a homogeneous linear equation with $k \ge 3;$ $c_1,\ldots, c_k \in \mathbb{Z} \setminus \{0\}$ such that $c_1+\ldots+c_k \neq 0,$ but there exists $\emptyset \neq S\subseteq [k]: \sum_{s\in S} c_s = 0$. Then, 
    $$
        d(\calL, \ve) = \Theta(\ve).
    $$
\end{conjecture}
\noindent We note that we would already find it interesting to determine whether the dependency on $k$ and/or $M$ is necessary or not in the exponent of $\ve$ here.

\noindent \textbf{Other limiting values.} Our classification theorem characterizes all homogeneous linear equations for which $d(\calL, \ve) \to 0$ as $\ve$ tends to zero. Following the classical Ramsey--Tur\'an theory, a very natural question to ask is what the limit of $d(\calL, \ve)$ is in case it is larger than $0$.

\begin{problem}
    Let $\calL$ be a homogeneous linear equation with no nonempty subset of coefficients with zero sum. Determine the value 
    $$
        \lim_{\ve \to 0} d(\calL, \ve)
    $$ if it exists.
\end{problem}

\noindent \textbf{Stronger assumptions on the independence number.} We note that in Ramsey--Tur\'an theory, there is another interesting direction concerned with quantitative bounds. Namely, what happens if we assume a stronger than just sublinear bound on the independence number? 
In our setting, this translates to asking, if we get more equations for which solution-free sets $A$ need to be sublinear if we assume, say $\alpha \le o(p/\log p)$? Or, if not, at what point do we start getting more such equations?

\noindent \textbf{Systems of equations.} There is a natural analogue of Roth's theorem for the density regularity of systems of equations. We note that this was considered already in Roth's original paper, but a full classification is a consequence of Szemer\'edi's theorem on arithmetic progressions (which can be thought of as solving a special case of the general density regularity problem for systems of equations). With this in mind, it is natural to ask what happens for our Ramsey--Tur\'an analogue when considering systems of linear equations.

\noindent \textbf{Other additive questions.} Finally, one can consider a similar Ramsey--Tura\'n variant of many other classical questions in additive combinatorics, and we suspect many lead to interesting directions for future research. In particular, since understanding for which problems examples which are close to tight are forced to be structured is a major area of modern additive combinatorics. We note here that using the notion of structure introduced by Erd\H{o}s and S\'ark\"ozy is perhaps the most natural analogue from the perspective of extending the Ramsey--Tur\'an theory from extremal combinatorics, but that other notions might lead to interesting questions as well.


\subsection*{Acknowledgement}
We want to thank Mehtaab Sawhney for useful discussions.
This research was conducted at Princeton University during the visits of the second, third, and fourth authors. They gratefully acknowledge Princeton University for its hospitality and support.

\providecommand{\MR}[1]{}
\providecommand{\MRhref}[2]{%
  \href{http://www.ams.org/mathscinet-getitem?mr=#1}{#2}
}

    \bibliographystyle{amsplain_initials_nobysame}
    \bibliography{bibfile}


\end{document}